\documentclass{conm-p-l}
\usepackage[colorlinks,
            linkcolor=reference,
            citecolor=citation,
            urlcolor=e-mail,
            backref]{hyperref}
\usepackage[all]{xy}
\SelectTips{cm}{}
\usepackage[mathscr]{euscript}
\usepackage{tikz}
\usepackage{mathabx}

\usepackage{color}
\definecolor{todo}{rgb}{1,0,0}
\definecolor{conditional}{rgb}{0,1,0}
\definecolor{e-mail}{rgb}{0,.40,.80}
\definecolor{reference}{rgb}{.20,.60,.22}
\definecolor{mrnumber}{rgb}{.80,.40,0}
\definecolor{citation}{rgb}{0,.40,.80}

\newcommand{\ie}{{\em i.e.}\ }
\newcommand{\cf}{{\em cf.}\ }

\newcommand{\ko}{\: , \;}

\newcommand{\ol}[1]{\overline{#1}}

\numberwithin{equation}{subsection}
\newtheorem{theorem}[subsection]{Theorem}
\newtheorem{classification-theorem}[subsection]{Classification Theorem}
\newtheorem{decomposition-theorem}[subsection]{Decomposition Theorem}
\newtheorem{proposition-definition}[subsection]{Proposition-Definition}
\newtheorem{periodicity-conjecture}[subsection]{Periodicity Conjecture}

\newtheorem{proposition}[subsection]{Proposition}

\newtheorem{example}[subsection]{Example}
\newtheorem{remark}[subsection]{Remark}

\newcommand{\reminder}[1]{}

\newcommand{\opname}[1]{\operatorname{\mathsf{#1}}}

\newcommand{\Mod}{\opname{Mod}\nolimits}

\newcommand{\tB}{\tilde{B}}

\newcommand{\cok}{\opname{cok}\nolimits}

\newcommand{\Z}{\mathbb{Z}}

\newcommand{\fg}{\mathfrak{g}}

\newcommand{\iso}{\xrightarrow{_\sim}}

\newcommand{\id}{\mathbf{1}}

\newcommand{\Hom}{\opname{Hom}}

\newcommand{\RHom}{\opname{RHom}}

\newcommand{\Ext}{\opname{Ext}}

\newcommand{\Tw}{\opname{Tw}}
\newcommand{\Alg}{\opname{Alg}}
\newcommand{\Coalg}{\opname{Coalg}}
\newcommand{\Set}{\opname{Set}}
\newcommand{\Com}{\opname{Com}}

\newcommand{\ten}{\otimes}

\newcommand{\cc}{{\mathscr C}}
\newcommand{\cd}{{\mathscr D}}

\newcommand{\eps}{\varepsilon}
\renewcommand{\phi}{\varphi}

\renewcommand{\tilde}[1]{\widetilde{#1}}

\makeatletter
\let\@wraptoccontribs\wraptoccontribs
\makeatother

\begin{document}

\date{October 23, 2019} 

\title{A remark on Hochschild cohomology and Koszul duality}

\author{Bernhard Keller}

\address{B. Keller: Universit\'e Paris Diderot -- Paris 7\\
    UFR de Math\'ematiques\\
    CNRS\\
   Institut de Math\'ematiques de Jussieu--Paris Rive Gauche, IMJ-PRG \\   
    B\^{a}timent Sophie Germain\\
    75205 Paris Cedex 13\\
    France
}
\email{bernhard.keller@imj-prg.fr}
\urladdr{https://webusers.imj-prg.fr/~bernhard.keller/}

\begin{abstract}
Applying recent results by Lowen--Van den Bergh we show that Hochschild
cohomology is preserved under Koszul--Moore duality as a Gerstenhaber
algebra. More precisely, the corresponding Hochschild complexes are
linked by a quasi-isomorphism of $B_\infty$-algebras.
\end{abstract}



\dedicatory{Dedicated to the Jos\'e Antonio de la Pe\~na on the occasion
of his sixtieth birthday}

\maketitle

\section{Introduction}
Consider the following statement:
\begin{quote}
Hochschild cohomology is preserved under Koszul duality.
\end{quote}
This is clearly wrong. Indeed, if $V$ is a non zero finite-dimensional vector space
over a field $k$, then the center (=zeroth Hochschild cohomology) of the 
symmetric algebra $SV$ is $SV$ but the center of the
Koszul dual exterior algebra $\Lambda(V^*)$ is finite-dimensional
(for a study of the Hochschild cohomology of $\Lambda(V^*)$, cf. \cite{Wong16}).
To try and save the statement, let us recall that $\Lambda(V^*)$ is in fact
the Yoneda algebra
\[
\Ext^*_{SV}(k,k)
\]
and that, with the zero differential,
it is even quasi-isomorphic to the derived endomorphism algebra
$\RHom_{SV}(k,k)$. Thus, it is natural to endow $\Lambda(V^*)$ with
the grading so that $V^*$ sits in degree $1$
and with the zero differential. With this new interpretation of $\Lambda(V^*)$
as a dg (=differential graded) algebra,
we find that its zeroth Hochschild cohomology is the {\em completion}
$\widehat{SV}$ of the symmetric algebra at the augmentation ideal.
It turns out that in order to get exactly $SV$, it suffices to replace
the dg algebra $\Lambda(V^*)$ with the $k$-dual dg {\em coalgebra} 
$\Lambda V$. In fact, we then get an algebra isomorphism
\[
HH^*(SV,SV) \iso HH^*(\Lambda V,\Lambda V).
\]
Our main results are that this isomorphism generalizes from $SV$ to
any augmented dg $k$-algebra (when we replace $\Lambda V$ with
the Koszul--Moore dual dg coalgebra) and that it lifts to an
isomorphism between the corresponding Hochschild cochain
complexes in the homotopy category of $B_\infty$-algebras.
In particular, the isomorphism in Hochschild cohomology is
an isomorphism of Gerstenhaber algebras. The algebra isomorphism
for augmented dg algebras follows from the setup of Koszul--Moore
duality, which we recall in section~\ref{s:Koszul-Moore duality}.
Our first proof of the lift to the $B_\infty$-level is now superseded
by recent work of Lowen--Van den Bergh \cite{LowenVandenBergh19},
which we use in the very short proof of Theorem~\ref{thm:B-infinity}.

Previous results relating the Hochschild cohomologies of Koszul(-Moore)
dual algebras can be found in \cite{Buchweitz03}, \cite{Keller03}
\cite{FelixMenichiThomas05}, \cite{BriggsGelinas17} and
\cite{BerglundBoerjeson17}
cf. remark~\ref{remark:history}.

\section*{Acknowledgments}
I am grateful to Jos\'e Antonio de la Pe\~{n}a for his friendship 
 throughout the years and for the
co-organization of many successful scientific events.
I thank the organizers of the ARTA 7, where this material was first presented.
I am indebted to Wendy Lowen for a talk at Trinity College Dublin
in May 2019 on the results of
\cite{LowenVandenBergh19}. Many thanks to Vladimir Dotsenko for reminding
me of references \cite{FelixMenichiThomas05} and \cite{BriggsGelinas17}
and to Pedro Tamaroff for pointing out \cite{BerglundBoerjeson17}.

\section{Reminder on Koszul--Moore duality}
\label{s:Koszul-Moore duality}

We follow Lef\`evre--Hasegawa \cite{Lefevre03} and Positselski
\cite{Positselski11}, \cf also \cite{Keller03a} and Appendix A to
\cite{VandenBergh15}.
Let $k$ be a field and $A$ a dg $k$-algebra. Thus,
$A$ is a $\Z$-graded associative algebra with $1$
\[
A=\bigoplus_{p\in\Z} A^p
\]
endowed with a homogeneous linear endomorphism $d$ of degree $1$, the
differential, such that $d^2=0$ and we have the Leibniz rule
\[
d(ab)=(da)b+(-1)^p a(db)
\]
for all $a\in A^p$ and all $b\in A$. Let $\eps: A \to k$ be an augmentation
(a morphism of dg algebras). For example, if $V$ is a vector space 
(concentrated in degree $0$), we can consider $A=SV$ (concentrated in
degree $0$ with $d=0$). Dually, let $C$ be a dg coalgebra and 
$\eps: k \to C$ a co-augmentation. For example, we may consider
$C=\Lambda V$, where $V$ is concentrated in degree $-1$ and $d=0$.
Denote by $\Hom_k(C,A)$ the graded vector space whose $n$th component
is formed by the homogeneous $k$-linear maps $f: C \to A$ of degree $n$.
We make $\Hom_k(C,A)$ into a differential graded algebra by setting
\[
d(f) = d\circ f - (-1)^n f\circ d
\]
for $f$ homogeneous of degree $n$ and
\[
f * g = \mu\circ (f\ten g)\circ \Delta
\]
for homogeneous $f$ and $g$, where $\mu$ is the multiplication of $A$ and
$\Delta$ the comultiplication of $C$. A {\em twisting cochain} is an element
$\tau\in \Hom_k(C,A)$ of degree $1$ such that
\[
\eps\circ\tau=0=\tau\circ \eta \mbox{ and } d(\tau)+\tau*\tau=0.
\]
For example, with
the above notations,  the composition of the natural projection and inclusion 
morphisms
\[
\Lambda V \to V \to SV
\]
is a twisting cochain. We denote by $\Tw(C,A)$ the set of twisting cochains. 
From now on, we assume that $C$ is {\em cocomplete}, \ie that
$\ol{C}=\cok(\eta)$ is the union of the kernels of the maps induced by
the iterated comultiplications
\[
\Delta^{(n)}: C \to C^{\ten n}\ko n\geq 2.
\]
We denote by $\Alg$ the category of augmented dg algebras and by
$\Coalg$ the category of cocomplete co-augmented dg coalgebras.

\begin{proposition}
\begin{itemize}
\item[a)] The functor $\Tw(?,A): \Coalg^{op} \to \Set$ is representable, \ie
there is an object $BA\in\Coalg$ and a functorial bijection
\[
\Tw(C',A) \iso \Coalg(C',BA).
\]
\item[b)] The functor $\Tw(C,?): \Alg \to \Set$ is co-representable, \ie there
is an object $\Omega A\in\Alg$ and a functorial bijection
\[
\Tw(C,A') \iso \Alg(\Omega C, A').
\]
\end{itemize}
\end{proposition}

The dg coalgebra $BA$ is known as the {\em bar construction} and the dg
algebra $\Omega C$ as the {\em cobar construction}. It is not hard
to describe $BA$ and $\Omega C$ explicitly but we will not need this.
Notice that if we
denote by $DBA$ the $k$-dual dg algebra of $BA$, then we have a
canonical isomorphism
\[
DBA \iso \RHom_A(k,k)
\]
and that the latter is quasi-isomorphic to the Koszul dual $A^!$ (with
the generators in differential degree $1$) if $A$ is a Koszul algebra
concentrated in degree $0$. We denote the category of dg right $A$-modules
by $\Mod A$. Its localization with respect to all quasi-isomorphisms is
the {\em derived category} $\cd A$. A (right) dg comodule $M$ is
{\em cocomplete} if it is the union of the kernels of the maps
\[
M \to M\ten \ol{C}^{n-1}\ko n\geq 2\ko
\]
induced by the iterated comultiplications. We denote by $\Com C$ the category
of cocomplete right dg $C$-comodules. It becomes a Frobenius exact category
when endowed with the conflations given by the exact sequences whose underlying
sequences of graded comodules split. The {\em category up to homotopy}
is the associated stable category. We denote by 
by $\cd C$ the {\em co-derived category}, \ie the localization of $\Com C$ 
at the class of all {\em co-quasi-isomorphisms}. Here, a morphism 
$s:L \to M$ of dg comodules
is a co-quasi-isomorphism if its cone lies in the smallest triangulated subcategory
of the category up to homotopy of dg comodules stable under coproducts and
containing all totalizations of short exact sequences
\[
0 \to X \to Y \to Z \to 0
\]
of dg comodules. This definition is due to Positselski \cite{Positselski11}.
One can also characterize the co-quasi-isomorphisms using the cobar
construction for dg comodules, \cf \cite{Lefevre03} and below. For the comparison
between the two, \cf Appendix~A of \cite{VandenBergh15}.

Let us fix a twisting cochain $\tau: C\to A$. For $M\in\Mod A$, let $M\ten_\tau C$
be the graded comodule $M\ten C$ endowed with the differential
\[
d\ten\id + \id\ten d +d_\tau \ko \mbox{ where } 
d_\tau=(\mu\ten\id)\circ (\id\ten\tau\ten\id) \circ (\id\ten\Delta).
\]
For $L\in\Com C$, let $L\ten_\tau A$ be the graded $A$-module $L\ten A$
endowed with the differential
\[
d\ten\id + \id\ten d +d_\tau \ko \mbox{ where } 
d_\tau=(\id\ten\mu)\circ (\id\ten\tau\ten\id) \circ (\Delta\ten\id).
\]

\begin{proposition} The pair $(?\ten_\tau A, ?\ten_\tau C)$ is a pair of
adjoint functors between $\Com C$ and $\Mod A$. It induces a pair of adjoint functors
between $\cd C$ and $\cd A$.
\end{proposition}

We define the twisting cochain $\tau: C \to A$ to be {\em acyclic} if the
associated adjoint functors between $\cd C$ and $\cd A$ are equivalences.
In this case, we say that $A$ and $C$ are {\em Koszul--Moore dual} to each other.

\begin{theorem}[\cite{Lefevre03}] \label{thm:lefevre} The following are equivalent:
\begin{itemize}
\item[i)] $\tau$ is acyclic;
\item[ii)] $\tau$ induces a quasi-isomorphism $\Omega C \to A$;
\item[iii)] $\tau$ induces a weak equivalence $C \to BA$ (\ie the induced
morphism $\Omega C \to \Omega BA$ is a quasi-isomorphism);
\item[iv)] The natural morphism $A\ten_\tau C \ten_\tau A \to A$ is a quasi-isomorphism.
\end{itemize}
\end{theorem}

For example, it follows from part iv) that $SV$ and $\Lambda V$ (with the above 
notations) are Koszul--Moore dual to each other. By part iii), we have Koszul--Moore
duality between $A$ and $BA$ and therefore a triangle equivalence
\[
\cd(BA) \iso \cd A.
\]
It is remarkable that the dg coalgebra $BA$ determines all of $\cd A$ whereas
the dual dg algebra $DBA \iso \RHom_A(k,k)$ a priori only determines the
thick subcategory of $\cd A$ generated by $k$.

Assume that $\tau: C \to A$ is an acylic cochain. Put $A^e=A\ten A^{op}$
and $C^e=C\ten C^{op}$ and let
\[
\tau^e=\tau\ten \eta + \eta \ten \tau: C^e \to A^e.
\]
Clearly $\tau^e$ is a twisting cochain.

\begin{proposition} \label{prop:alg-isom}
The twisting cochain $\tau^e$ is acyclic and the
induced equivalence 
\[
\cd(C^e) \iso \cd(A^e)
\]
takes the dg bicomodule $C$ to the dg bimodule $A$. Thus we have an induced
isomorphism of graded algebras
\[
HH^*(C)=\Ext^*_{C^e}(C,C) \iso \Ext^*_{A^e}(A,A) = HH^*(A).
\]
\end{proposition}

\begin{remark}
We see that Koszul--Moore duality preserves Hochschild cohomology as 
a graded algebra.
\end{remark}

\begin{proof} Let $\phi_A : A \ten_\tau C \ten_\tau A \to A$ be the canonical
morphism. Then the canonical morphism
\[
\phi_{A^e} : A^e \ten_{\tau^e} C^e \ten_{\tau^e} A^e \to A^e
\]
is the composition of the isomorphism
\[
(A\ten A^{op})\ten_{\tau^e} (C\ten C^{op}) \ten_{\tau^e} (A\ten A^{op}) \iso
(A\ten_\tau C \ten_\tau A) \ten (A^{op} \ten_\tau C^{op} \ten_\tau A^{op})
\]
with the quasi-isomorphism
\[
\phi_A \ten \phi_{A^{op}} : (A\ten_\tau C \ten_\tau A) \ten (A^{op} \ten_\tau C^{op} \ten_\tau A^{op}) \to A\ten A^{op}.
\]
Thus $\tau^e$ is an ayclic twisting cochain by Theorem~\ref{thm:lefevre}. 
The induced equivalence takes $C$ to 
\[
C\ten_{\tau^e}(A\ten A^{op}) = A \ten_\tau C \ten_\tau A
\]
which is quasi-isomorphic to $A$ since $\tau$ is acyclic.
\end{proof}

\begin{theorem} \label{thm:B-infinity}
The isomorphism of the proposition lifts to an isomorphism
in the homotopy category of $B_\infty$-algebras between the corresponding
Hochschild cochain complexes. In particular, it preserves the Gerstenhaber
brackets.
\end{theorem}

\begin{remark} \label{remark:history}
In \cite{Buchweitz03}, Buchweitz related the Hochschild cohomology algebras of
a Koszul algebra and its Koszul dual.
In Theorem~3.5 of \cite{Keller03}, we showed that for a 
{\em Koszul algebra} $A$ (concentrated in differential degree $0$), there is a 
canonical isomorphism in the category of Adams graded $B_\infty$-algebras
between the Hochschild complex of $A$ and that of the Koszul dual
algebra $A^!$ whose $p$th graded piece is put into bidegree $(p,-p)$, 
where the first component is the differential degree and the second
component the Adams degree. Here, we get rid of the Koszulity
assumption and the Adams grading by using dg coalgebras.

In \cite{FelixMenichiThomas05}, the F\'elix--Menichi--Thomas show that for a simply
connected coalgebra $C$, the Hochschild cohomologies of the
dg algebras $\Omega C$ and $\Hom_k(C,k)$ are isomorphic as
Gerstenhaber algebras. Another proof of this, under less stringent
connectedness assumptions, is given in section~3 of 
Briggs--G\'elinas' \cite{BriggsGelinas17}. For Koszul $A_\infty$--algebras,
an isomorphism of the Hochschild cohomologies as weight graded
$A_\infty$-algebras (but not as Gerstenhaber algebras) is proved
by Berglund--B\"orjeson in Theorem~3.2 of \cite{BerglundBoerjeson17}.
\end{remark} 

\begin{proof} The twisting cochain induces a weak equivalence
$BA \to C$. By the argument of \cite{Keller03}, this yields an isomorphism 
in the homotopy category of
$B_\infty$-algebras between the Hochschild complexes of $BA$ and $C$
(these are sometimes called coHochschild complexes).
Thus, we may assume that $C=BA$ and $\tau$ is the canonical twisting
cochain. Let $\tB A=A\ten_\tau C \ten_\tau A$. Then $\tB A$ has a natural
structure of dg coalgebra in the category $\Mod(A^e)$ of dg $A$-$A$-bimodules
endowed with $\ten_A$. We have a lax monoidal functor
\[
F: \Mod(A^e) \to \cc(k)
\]
taking a bimodule $M$ to $k\ten_A M\ten_A k$. The lax structure is given by
the morphism
\[
F(L\ten_A M) =k\ten_A (L\ten_A M)\ten_A k = k\ten_A (L \ten_A A \ten_A M) \ten_A k \to
k\ten_A L \ten_A k \ten_A M\ten_A k = (FL)\ten_k (FM).
\]
The functor $F$ sends the $A^e$-coalgebra $\tB A$ to the $k$-coalgebra $BA=C$.
It induces a morphism from the Hochschild complex of $\tB A$ to that
of $C$ and this is easily seen to be a quasi-isomorphism compatible with
the cup product and the brace operations. The claim follows by Theorem~5.1
of \cite{LowenVandenBergh19} which states that there is a canonical
isomorphism in the homotopy category of $B_\infty$-algebras between
the Hochschild complex of $\tB A$ and the Hochschild complex of $A$.
It is not hard to check that in homology, it induces the isomorphism
of Propostion~\ref{prop:alg-isom}.
\end{proof}

\begin{example} Suppose that $k$ is of characteristic $0$ and $\fg$ a finite-dimensional
Lie algebra over $k$. Let $U\fg$ be its enveloping algebra and $\Lambda\fg$
the supersymmetric coalgebra on $\fg$ placed in degree $-1$ and endowed with
the coalgebra differential whose $(-2)$-component is the bracket $\Lambda^2\fg \to \fg$.
Thus, the underlying complex of $\Lambda \fg$ is the Chevalley--Eilenberg complex
of $\fg$.
Then the map $\tau: \Lambda\fg \to U\fg$ which is the composition of the
projection $\Lambda\fg \to \fg$ with the inclusion $\fg \to U\fg$ is an acyclic
twisting cochain and we obtain an isomorphism of Gerstenhaber algebras
\[
HH^*(\Lambda \fg) \iso HH^*(U\fg).
\]
It would be interesting to generalize this example to Lie algebroids.
\end{example}



\def\cprime{$'$} \def\cprime{$'$}
\providecommand{\bysame}{\leavevmode\hbox to3em{\hrulefill}\thinspace}
\providecommand{\MR}{\relax\ifhmode\unskip\space\fi MR }
\providecommand{\MRhref}[2]{%
  \href{http://www.ams.org/mathscinet-getitem?mr=#1}{#2}
}
\providecommand{\href}[2]{#2}

\end{document}